\documentclass[natbib]{svcyclop}
\bibpunct{[}{]}{;}{n}{}{,} 
\usepackage{latexsym}
\usepackage{graphicx}

\usepackage{amssymb}
\usepackage{amsmath}

\newcommand{\mb}[1]{\mathbf{#1}}

\usepackage{ulem}

\begin{document}

\title{Derivative based global sensitivity measures}

\author{Serge\"i Kucherenko \inst{1} and Bertrand Iooss \inst{2,3} }

\institute{
Imperial College London\\
 London, SW7 2AZ, UK\\
E-mail: s.kucherenko@imperial.ac.uk
\and
EDF R\&D\\
6 quai Watier, 78401 Chatou, France\\
E-mail: bertrand.iooss@edf.fr
\and
Institut de Math\'ematiques de Toulouse\\
Universit\'e Paul Sabatier\\
118 route de Narbonne, 31062 Toulouse, France}

\maketitle

\section{Abstract}

The method of derivative based global sensitivity measures (DGSM) has recently become popular among practitioners. It has a strong link with the Morris screening method and Sobol' sensitivity indices and has several advantages over them. DGSM are very easy to implement and evaluate numerically. The computational time required for numerical evaluation of DGSM is generally much lower than that for estimation of Sobol' sensitivity indices. This paper presents a survey of recent advances in DGSM concerning lower and upper bounds on the values of Sobol' total sensitivity indices $S_{i}^{\mbox{\scriptsize{tot}}}$. Using these bounds it is possible in most cases to get a good practical estimation of the values of $S_{i}^{\mbox{\scriptsize{tot}}} $. Several examples are used to illustrate an application of DGSM.

\vspace{0.5cm}
Keywords: Sensitivity analysis, Sobol' indices, Morris method, Model derivatives, DGSM, Poincar\'e inequality

\section{Introduction}

Global sensitivity analysis (SA) offers a comprehensive approach to the model analysis. Unlike local SA, global SA methods evaluate the effect of a factor while all other factors are varied as well and thus they account for interactions between variables and do not depend on the choice of a nominal point. Reviews of different global SA methods can be found in \citet{salrat08} and \citet{sobkuc05}. The method of global sensitivity indices suggested by \citet{sob90,sob93}, and then further developed by \citet{homsal96} is one of the most efficient and popular global SA techniques. It belongs to the class of variance-based methods. These methods provide information on the importance of different subsets of input variables to the output variance. There are two types of Sobol' sensitivity indices: the main effect indices, which estimate the individual contribution of each input parameter to the output variance, and the total sensitivity indices, which measure the total contribution of a single input factor or a group of inputs. The total sensitivity indices are used to identify non-important variables which can then be fixed at their nominal values to reduce model complexity. This approach is known as ``factors' fixing setting'' \citep{salrat08}. For high-dimensional models the direct application of variance-based global SA measures can be extremely time-consuming and impractical.

A number of alternative SA techniques have been proposed. One of them is the screening method by \citet{mor91}. It can be regarded as global as the final measure is obtained by averaging local measures (the elementary effects). This method is considerably cheaper than the variance based methods in terms of computational time. The Morris method can be used for identifying unimportant variables. However, the Morris method has two main drawbacks. Firstly, it uses random sampling of points from the fixed grid (levels) for averaging elementary effects which are calculated as finite differences with the increment delta comparable with the range of uncertainty. For this reason it can not correctly account for the effects with characteristic dimensions much less than delta. Secondly, it lacks the ability of the Sobol' method to provide information about main effects (contribution of individual variables to uncertainty) and it can't distinguish between low and high order interactions.

This paper presents a survey of derivative based global sensitivity measures (DGSM) and their link with Sobol' sensitivity indices. DGSM are based on averaging local derivatives using Monte Carlo or Quasi Monte Carlo sampling methods. This technique is much more accurate than the Morris method as the elementary effects are evaluated as strict local derivatives with small increments compared to the variable uncertainty ranges. Local derivatives are evaluated at randomly or quasi randomly selected points in the whole range of uncertainty and not at the points from a fixed grid. 

The so-called alternative global sensitivity estimator defined as a normalized integral of partial derivatives was firstly introduced by \citet{sobger95}.
\citet{kucrod09} introduced some other DGSM and coined the acronym DGSM. They showed that DGSM can be seen as the generalization of the Morris method \citep{mor91}. 
\citet{kucrod09}  also established empirically the link between DGSM and Sobol' sensitivity indices. They showed that the computational cost of numerical evaluation of DGSM can be much lower than that for estimation of Sobol' sensitivity indices. 

\citet{sobkuc09}  proved theoretically that, in the cases of uniformly and normally distributed input variables, there is a link between DGSM and the Sobol' total sensitivity index $S_{i}^{\mbox{\scriptsize{tot}}} $ for the same input. They showed that DGSM can be used as an upper bound on total sensitivity index $S_{i}^{\mbox{\scriptsize{tot}}} $. Small values of DGSM  imply small $S_{i}^{\mbox{\scriptsize{tot}}} $, and hence unessential factors $x_{i} $. 
However, ranking influential factors using DGSM can be similar to that based on $S_{i}^{\mbox{\scriptsize{tot}}} $ only for the case of linear and quasi-linear models. For highly non-linear models two rankings can be very different. They also introduced modified DGSM which can be used for both a single input and groups of inputs \citep{sobkuc10}. 
From DGSM, \citet{kucson15} have also derived lower bounds on total sensitivity index.
\citet{lamioo13} extended results of Sobol' and Kucherenko for models with input variables belonging to the general class of continuous probability distributions.
In the same framework, \citet{roufru14} have defined crossed-DGSM, based on second-order derivatives of model output, in order to bound the total Sobol' indices of an interaction between two inputs.

All these DGSM measures can be applied for problems with a high number of input variables to reduce the computational time. 
Indeed, the numerical efficiency of the DGSM method can be improved by using the automatic differentiation algorithm for calculation DGSM as was shown in \citet{kipkuc09}. 
However, the number of required function evaluations still remains to be proportional to the number of inputs. 
This dependence can be greatly reduced using an approach based on algorithmic differentiation in the adjoint or reverse mode \citep{griwal08} (      \uline{Variational Methods}). 
It allows estimating all derivatives at a cost at most 4-6 times of that for evaluating the original function \citep{janleo14}.

This paper is organised as follows: the Morris method and DGSM are firstly described in the following section. 
Sobol' global sensitivity indices and useful relationships are then introduced.
Therefore, DGSM-based lower and uppers bounds on total Sobol' sensitivity indices for uniformly and normally distributed random variables are presented, followed by DGSM for groups of variables and their link with total Sobol' sensitivity indices.
Another section presents the upper bounds results in the general case of variables with continuous probability distributions.
Then, computational costs are considered, followed by some test cases which illustrate an application of DGSM and their links with total Sobol' sensitivity indices. 
Finally, conclusions are presented in the last section.

\section{From Morris method to DGSM}
 
\subsection{Basics of the Morris method}

The Morris method is traditionally used as a screening method for problems with a high number of variables for which function evaluations can be CPU-time consuming (see \uline{Design of Experiments for Screening}). It is composed of individually randomized 'one-factor-at-a-time' (OAT) experiments. 
Each input factor may assume a discrete number of values, called levels, which are chosen within the factor range of variation.

The sensitivity measures proposed in the original work of \citet{mor91} are based on what is called an elementary effect. It is defined as follows. The range of each input variable is divided into \textit{p} levels. Then the elementary effect (incremental ratio) of the \textit{i-}th input factor is defined as
\begin{equation} \label{GrindEQ__2_} 
EE_{i} \left(\mb{x}^{*} \right)=\frac{\left[G\left(x_{1}^{*} ,\ldots ,x_{i-1}^{*} ,x_{i}^{*} +\Delta ,x_{i+1}^{*} ,\ldots ,x_{d}^{*} \right)-G\left(\mb{x}^{*} \right)\right]}{\Delta } ,  
\end{equation} 
where $\Delta $ is a predetermined multiple of 1/(\textit{p}-1) and point $\mb{x}^{*}=(x_1^{*},\ldots,x_d^{*}) \in H^{d} {\rm \; }$ is such that $x_{i}^{*} +\Delta \le 1$. One can see that the elementary effect are finite difference approximations of the model derivative with respect to $x_i$ and using a large perturbation step $\Delta$.

The distribution of elementary effects $EE_i$ is obtained by randomly sampling \textit{R} points from $H^{d} $. Two sensitivity measures are evaluated for each factor: $\mu_{i} $ an estimate of the mean of the distribution $EE_i$, and $\sigma_{i} $ an estimate of the standard deviation of $EE_i$. A high value of $\mu_{i} $ indicates an input variable with an important overall influence on the output. A high value of $\sigma_{i} $ indicates a factor involved in interaction with other factors or whose effect is nonlinear. The computational cost of the Morris method is \textit{N${}_{F }$} = \textit{R} (\textit{d+1}). 

The revised version of the $EE_{i} \left(\mb{x}^{*} \right)$ measure and a more effective sampling strategy, which allows a better exploration of the space of the uncertain input factors was proposed by \citet{camcar07}. 
To avoid the canceling effect which appears in non-monotonic functions \citet{camcar07} introduced another sensitivity measure $\mu_{i}^{*} $ based on the absolute value of $EE_{i} (\mb{x}^{*})$: $\left|EE_{i} (\mb{x}^{*})\right|$. 
It was also noticed that $\mu_{i}^{*} $ has similarities with the total sensitivity index $S_{i}^{\mbox{\scriptsize{tot}}} $ in that it can give a ranking of the variables similar to that based on the $S_{i}^{\mbox{\scriptsize{tot}}} $ but no formal proof of the link between $\mu_{i}^{*} $ and $S_{i}^{\mbox{\scriptsize{tot}}}$ was given \citep{camcar07}. 

%

Finally, other extensions of the initial Morris method have been introduced for the second-order effects' analysis \citep{cambra99} \citep{crobra02} \citep{fedran15}, for the estimation of Morris' measures with any-type of design \citep{puj09} \citep{sancor12} and for building some 3D Morris' graph \citep{puj09}.

\subsection{The local sensitivity measure}
 
Consider a differentiable function $G\left(\mb{x}\right)$, where $\mb{x}=(x_1,\ldots,x_d)$ is a vector of input variables defined in the unit hypercube $H^{d} $ $\left(0\le x_{i} \le 1\, ,\; i=1,\ldots ,d\right)$. 
Local sensitivity measures are based on partial derivatives 
\begin{equation} \label{GrindEQ__1_} 
E_{i} (\mb{x}^{*})=\frac{\partial G(\mb{x}^{*})}{\partial x_{i} } .      
\end{equation} 
This measure $E_i$ is the limit version of the elementary effect $EE_i$ defined in \eqref{GrindEQ__1_} when $\Delta$ tends to zero.
It is its generalization in this sense.
In SA, using the partial derivative $\partial G \left/ \partial x_{i}\right.$ is well known as a local method (see \uline{Variational Methods}).
In this paper, the goal is to take advantage of this information in global SA.

The local sensitivity measure $E_{i} (\mb{x}^{*})$ depends on a nominal point $\mb{x}^*$ and it changes with a change of $\mb{x}^{*}$. 
This deficiency can be overcome by averaging $E_{i} (\mb{x}^{*})$ over the parameter space $H^{d}$.
This is done just below, allowing to define new sensitivity measures, called DGSM for Derivative-based Global Sensitivity Measures.

\subsection{DGSM for uniformly distributed variables}

Assume that ${\partial G\mathord{\left/ {\vphantom {\partial G \partial x_{i} }} \right. \kern-\nulldelimiterspace} \partial x_{i} } \in L_{2} $. 
Three different DGSM measures are defined:
\begin{equation} \label{GrindEQ__10_} 
\nu_{i} =\int_{H^{d} }\left(\frac{\partial G(\mb{x})}{\partial x_{i} } \right)^{2} d\mb{x} ,
\end{equation} 
\begin{equation} \label{GrindEQ__11_} 
w_{i}^{(m)} =\int_{H^{d} }x_{i}^{m} \frac{\partial G(\mb{x})}{\partial x_{i} } d \mb{x},     
\end{equation} 
where $m>0$ is a constant, and
\begin{equation} \label{GrindEQ__12_} 
\varsigma_{i} =\frac{1}{2} \int_{H^{d} }x_{i} (1-x_{i} )\left(\frac{\partial G(\mb{x})}{\partial x_{i} } \right)^{2} d \mb{x}.    
\end{equation} 

\subsection{DGSM for randomly distributed variables}

Consider a function $G\left(X_{1} ,...,X_{d} \right)$, where $X_{1} ,...,X_{d} $ are independent random variables, defined in the Euclidian space $R^{d}$, with cumulative density functions (cdfs) $F_{1} \left(x_{1} \right),...,F_{d} \left(x_{d} \right)$. 
The following DGSM was introduced in \citet{sobkuc09}:
\begin{equation} \label{GrindEQ__36_} 
\nu_{i} =\int_{R^{d} }\left(\frac{\partial G(\mb{x})}{\partial x_{i} } \right)^{2} d F(\mb{x}) = \mathbb{E}\left[\left(\frac{\partial G(\mb{x})}{\partial x_{i} } \right)^{2}\right],
\end{equation} 
with $F$ the joint cdf. A new measure is also introduced:
\begin{equation} \label{GrindEQ__37_} 
w_{i} =\int_{R^{d} }\frac{\partial G(\mb{x})}{\partial x_{i} } d F(\mb{x})  = \mathbb{E}\left(\frac{\partial G(\mb{x})}{\partial x_{i} } \right).
\end{equation} 

In \eqref{GrindEQ__10_} and \eqref{GrindEQ__36_}, $\nu_{i} $ is in fact the mean value of $\left({\partial G\mathord{\left/ {\vphantom {\partial G \partial x_{i} }} \right. \kern-\nulldelimiterspace} \partial x_{i} } \right)^{2} $.
In the following and in practice, it will be the most useful DGSM.

\section{Sobol' global sensitivity indices}

\subsection{Definitions}

The method of global sensitivity indices developed by Sobol' (see \uline{Variance-based Sensitivity Analysis: Theory and Estimation Algorithms}) is based on ANOVA decomposition \citep{harlit73}. 
Consider a square integrable function $G(\mb{x})$ defined in the unit hypercube $H^{d} $. It can be expanded in the following form
\begin{equation} \label{GrindEQ__3_} 
G(\mb{x})=g_{0} +\sum_{i}g_{i} (x_{i} )+\sum_{i<j}g_{ij} (x_{i} ,x_{j} )+...+g_{12...d} (x_{1} ,x_{2} ,...,x_{d} )  .   
\end{equation} 
This decomposition is unique if conditions $\displaystyle \int_{0}^{1} g_{i_{1} ...i_{s} } dx_{i_{k} }  =0 \quad$ for  ${\rm 1}\le k\le s$,  are satisfied. Here $1 \leq \textit{i}_{1} <\dots < \textit{i}_{s} \leq d$.

The variances of the terms in the ANOVA decomposition add up to the total variance of the function
\[V=\sum_{s=1}^{d}\sum_{i_{1} <\cdot \cdot \cdot <i_{s} }^{d}V_{i_{1} ...i_{s} }   ,    \] 
where $\displaystyle V_{i_{1} ...i_{s} } =\int_{0}^{1}g_{i_{1} ...i_{s} }^{2} (x_{i_{1} } ,...,x_{i_{s} } )dx_{i_{1} } ,...,x_{i_{s} }  $ are called partial variances.

Sobol' defined the global sensitivity indices as the ratios
\[S_{i_{1} ...i_{s} } =V_{i_{1} ...i_{s} } /V.\] 
All $S_{i_{1} ...i_{s} } $ are non negative and add up to one:
\[\sum_{i=1}^{d}S_{i}  +\sum_{i}\sum_{j}S_{ij}   +\sum_{i}\sum_{j}\sum_{k}S_{ijk}    ...+S_{1,2,...,d} =1 .\] 
Sobol' also defined sensitivity indices for subsets of variables. Consider two complementary subsets of variables $y$ and $z$:
\[\mb{x}=(y,z).     \] 
Let $y=(x_{i_{1} } ,...,x_{i_{m} } ), 1\le i_{1} <...<i_{m} \le d, K=(i_{1} ,...,i_{m} )$. The variance corresponding to the set $y$ is defined as
\[V_y =\sum_{s=1}^{m}\sum_{(i_{1} <\cdot \cdot \cdot <i_{s} )\in K}V_{i_{1} ...i_{s} }   .     \] 
$V_{y} $ includes all partial variances $V_{i_{1} } $, $V_{i_{2} } $,\dots , $V_{i_{1} ...i_{s} } $ such that their subsets of indices $(i_{1} ,...,i_{s} )\in K$. 

The total sensitivity indices were introduced by \citet{homsal96}. The total variance $V_{y}^{\mbox{\scriptsize{tot}}} $ is defined as 
\[V_{y}^{\mbox{\scriptsize{tot}}} =V-V_{z} .      \] 
$V_{y}^{\mbox{\scriptsize{tot}}} $ consists of all $V_{i_{1} ...i_{s} } $ such that at least one index $i_{p} \in K$ while the remaining indices can belong to the complimentary to \textit{K} set $\bar{K}$. The corresponding global sensitivity indices are defined as
\begin{equation} \label{GrindEQ__4_} 
\begin{array}{l} {S_{y} =V_{y} /V,} \\ {S_{y}^{\mbox{\scriptsize{tot}}} =V_{y}^{\mbox{\scriptsize{tot}}} /V.} \end{array} 
\end{equation} 

The important indices in practice are $S_{i} $ and $S_{i}^{\mbox{\scriptsize{tot}}}$, $i=1,...,d$:
\begin{equation} \label{GrindEQ__5_} 
\begin{array}{l} {S_{i} =V_{i} /V,} \\ {S_{i}^{\mbox{\scriptsize{tot}}} =V_{i}^{\mbox{\scriptsize{tot}}} /V.} \end{array} 
\end{equation} 
Their values in most cases provide sufficient information to determine the sensitivity of the analyzed function to individual input variables. Variance-based methods generally require a large number of function evaluations (see \uline{Variance-based Methods: Theory and Algorithms}) to achieve reasonable convergence and can become impractical for large engineering problems.

\subsection{Useful relationships}

To present further results on lower and upper bounds of $S_{i}^{\mbox{\scriptsize{tot}}}$, new notations and useful relationships have to be firstly presented. 
Denote $u_{i}(\mb{x})$ the sum of all terms in the ANOVA decomposition \eqref{GrindEQ__3_} that depend on $x_{i} $:
\begin{equation} \label{GrindEQ__13_} 
u_{i} (\mb{x})=g_{i} (x_{i} )+\sum_{j=1,j\ne i}^{d}g_{ij} (x_{i} ,x_{j} ) +\cdots +g_{12\cdots d} (x_{1} ,\cdots ,x_{d} ).   
\end{equation} 
From the definition of ANOVA decomposition it follows that
\begin{equation} \label{GrindEQ__14_} 
\int_{H^{d} }u_{i} (\mb{x})d\mb{x} =0.       
\end{equation} 
It is obvious that 
\begin{equation} \label{GrindEQ__15_} 
\frac{\partial G}{\partial x_{i} } =\frac{\partial u_i}{\partial x_{i} } .
\end{equation} 
Denote $\mb{z}=(x_{1} ,...,x_{i-1} ,x_{i+1} ,...,x_{d} )$ the vector of all variables but $x_{i} $, then $\mb{x} \equiv (x_{i} ,\mb{z})$ and $G(\mb{x}) \equiv G(x_{i} ,\mb{z})$. The ANOVA decomposition of $G(\mb{x})$ \eqref{GrindEQ__3_} can be presented in the following form
\[G(\mb{x})=u_{i} (x_{i} ,\mb{z})+v(\mb{z}),      \] 
where $v(\mb{z})$ is the sum of terms independent of $x_{i} $. Because of \eqref{GrindEQ__14_} it is easy to show that $\displaystyle v(\mb{z})=\int_0^1 G(\mb{x})dx_{i}$. 
Hence
\begin{equation} \label{GrindEQ__16_} 
u_{i} (x_{i} ,\mb{z})=G(\mb{x})-\int_0^1 G(\mb{x})dx_{i} .
\end{equation} 
This equation can be found in \citet{lam13}.
The total partial variance $V_{i}^{\mbox{\scriptsize{tot}}} $ can be computed as 
\[V_{i}^{\mbox{\scriptsize{tot}}} =\int_{H^{d} }u_{i}^{2} (\mb{x})d\mb{x} =\int_{H^{d} }u_{i}^{2} (x_{i} ,z)dx_{i} d\mb{z} .   \] 
Then the total sensitivity index $S_{i}^{\mbox{\scriptsize{tot}}}$ \eqref{GrindEQ__5_} is equal to
\begin{equation} \label{GrindEQ__17_} 
S_{i}^{\mbox{\scriptsize{tot}}} =\frac{1}{V} \int_{H^{d} }u_{i}^{2} (\mb{x})d\mb{x} .      
\end{equation} 

\subsection{A first direct link between total Sobol' sensitivity indices and partial derivatives}

Consider continuously differentiable function $G(\mb{x})$ defined in the unit hypercube $H^{d} $=$[0,1]^{d} $. 
This section presents a theorem that establishes links between the index $S_{i}^{\mbox{\scriptsize{tot}}} $ and the limiting values of $\left|{\partial G\mathord{\left/ {\vphantom {\partial G \partial x_{i} }} \right. \kern-\nulldelimiterspace} \partial x_{i} } \right|$.

In the case when $\mb{y}=\left(x_{i} \right)$, Sobol'-Jansen formula \citep{jan99}\citep{sob01}\citep{salann09} for $D_{i}^{\mbox{\scriptsize{tot}}} $ can be rewritten as
\begin{equation} \label{GrindEQ__6_} 
D_{i}^{\mbox{\scriptsize{tot}}} =\frac{1}{2} \int_{H^{d} }\int_{0}^{1}\left[G\left(\mb{x}\right)-G\left(\mathop{\mb{x}}\limits^{\circ } \right)\right]^{2} d\mb{x}dx'_{i}   ,    
\end{equation} 
where $\mathop{\mb{x}}\limits^{{\it o}} =\left(x_{1} ,...,x_{i-1} ,x'_{i} ,x_{i+1} ,...,x_{n} \right)$.    

\noindent \textbf{Theorem 1.} Assume that $\displaystyle c\le \left|\frac{\partial G}{\partial x_{i} } \right|\le C$, then
\begin{equation} \label{GrindEQ__7_} 
\frac{c^{2} }{12V} \le S_{i}^{\mbox{\scriptsize{tot}}} \le \frac{C^{2} }{12V} .      
\end{equation} 

\noindent \textbf{Proof:} Consider the increment of $G\left(\mb{x}\right)$ in \eqref{GrindEQ__6_}:
\begin{equation} \label{GrindEQ__8_} 
G\left(\mb{x}\right)-G\left(\mathop{\mb{x}}\limits^{\circ } \right)=\frac{\partial G\left(\hat{\mb{x}}\right)}{\partial x_{i} } \left(x_{i} -x'_{i} \right),     
\end{equation} 
where $\hat{\mb{x}}$ is a point between $\mb{x}$ and $\mathop{\mb{x}}\limits^{\circ } $. Substituting \eqref{GrindEQ__8_} into \eqref{GrindEQ__6_} leads to
\begin{equation} \label{GrindEQ__9_} 
V_{i}^{\mbox{\scriptsize{tot}}} =\frac{1}{2} \int_{H^{d} }\int_{0}^{1}\left(\frac{\partial G\left(\hat{\mb{x}}\right)}{\partial x_{i} } \right)^{2} \left(x_{i} -x'_{i} \right)^{2} d\mb{x}dx'_{i}   .    
\end{equation} 
In \eqref{GrindEQ__9_} $c^{2} \le \left({\partial G\mathord{\left/ {\vphantom {\partial G \partial x_{i} }} \right. \kern-\nulldelimiterspace} \partial x_{i} } \right)^{2} \le C^{2} $ while the remaining integral is 
\[\int_{0}^{1}\int_{0}^{1}\left(x'_{i} -x_{i} \right)^{2} dx'_{i} dx_{i}  =\frac{1}{6}  .    \] 
Thus obtained inequalities are equivalent to \eqref{GrindEQ__7_}. Consider the function $G=g_{0} +c(x_i-1/2)$. In this case $C=c$, $V=1/12$ and $S_{i}^{\mbox{\scriptsize{tot}}} =1$ and the inequalities in \eqref{GrindEQ__7_} become equalities.
{\hfill\mbox{\rule{2 true mm}{3 true mm}}}

\section{DGSM-based bounds for uniformly and normally distributed variables}

In this section, several theorems are listed in order to define useful lower and upper bounds of the total Sobol' indices.
The proofs of these theorems come from previous works and papers and are not recalled here.
Two cases are considered: variables $\mb{x}$ following uniform distributions and variables $\mb{x}$ following Gaussian distributions.
The general case will be seen in a subsequent section.

\subsection{Uniformly distributed variables}

\subsubsection{Lower bounds on $S_{i}^{\mbox{\scriptsize{tot}}} $}

\noindent \textbf{Theorem 2. }There exists the following lower bound between DGSM \eqref{GrindEQ__10_} and the Sobol' total sensitivity index:
\begin{equation} \label{GrindEQ__18_} 
\frac{\left(\int_{H^{d} }\left[G\left(1,\mb{z}\right)-G\left(0,\mb{z}\right)\right]\left[G\left(1,\mb{z}\right)+G\left(0,\mb{z}\right)-2G\left(\mb{x}\right)\right]d \mb{x}\right)^{2} }{4\nu_{i} V} <S_{i}^{\mbox{\scriptsize{tot}}}  
\end{equation} 

\noindent \textbf{Proof:} 
The proof of this Theorem is given in \citet{kucson15} and is based on equation \eqref{GrindEQ__17_} and a Cauchy-Schwartz inequality applied on $\displaystyle \int_{H^{d} }u_{i} (\mb{x})\frac{\partial u_{i} (\mb{x})}{\partial x_{i} }  d\mb{x}$.
{\hfill\mbox{\rule{2 true mm}{3 true mm}}}

The lower bound number number one (LB1) is defined as
\[\frac{\left(\int_{H^{d} }\left[G\left(1,\mb{z}\right)-G\left(0,\mb{z}\right)\right]\left[G\left(1,\mb{z}\right)+G\left(0,\mb{z}\right)-2G\left(\mb{x}\right)\right]d \mb{x}\right)^{2} }{4\nu_{i} V} .\]

\noindent \textbf{Theorem 3. }There exists the following lower bound, denoted $\gamma(m)$, between DGSM \eqref{GrindEQ__11_} and the Sobol' total sensitivity index:
\begin{equation} \label{GrindEQ__23_} 
\gamma(m) = \frac{(2m+1)\left[\int_{H^{d} }\left(G(1,\mb{z})-G(\mb{x})\right)d\mb{x} -w_{i}^{(m+1)} \right]^{2} }{(m+1)^{2} V} <S_{i}^{\mbox{\scriptsize{tot}}} .     
\end{equation} 

\noindent \textbf{Proof:} 
The proof of this Theorem in given in \citet{kucson15}  and is based on equation \eqref{GrindEQ__17_} and a Cauchy-Schwartz inequality applied on $\displaystyle \int_{H^{d} }x_{i}^{m} u_{i} (\mb{x})d \mb{x}$.
{\hfill\mbox{\rule{2 true mm}{3 true mm}}}
In fact, Theorem 3 gives a set of lower bounds depending on parameter \textit{m}. 
The value of \textit{m} at which $\gamma (m)$ attains its maximum is of particular interest. Further, star ($^*$) is used to denote such a value $m$: $m^{*} =\arg  \max (\gamma (m))$.
$\gamma(m^*)$ is called the lower bound number two (LB2):
\begin{equation} \label{GrindEQ__29_} 
\gamma(m^{*} )=\frac{(2m^{*} +1)\left[\int_{H^{d} }\left(G(1,\mb{z})-G(\mb{x})\right)d\mb{x} -w_{i}^{(m^{*} +1)} \right]^{2} }{(m^{*} +1)^{2} V}  
\end{equation} 

The maximum lower bound LB* is defined as 
\begin{equation}\label{GrindEQ__30_}
\mbox{LB*}=\max(\mbox{LB1,LB2}).    
\end{equation}  
Both lower and upper bounds can be estimated by a set of derivative based measures:
\begin{equation} \label{GrindEQ__31_} 
\Upsilon_{i} =\{ \nu_{i},w_{i}^{(m)},\zeta_i \} ,{\rm \; }m>0.       
\end{equation} 

\subsubsection{Upper bounds on $S_{i}^{\mbox{\scriptsize{tot}}} $ }

\noindent \textbf{\underbar{}}

\noindent \textbf{Theorem 4.} There exists the following upper bound between DGSM \eqref{GrindEQ__10_} and the Sobol' total sensitivity index:
\begin{equation} \label{GrindEQ__32_} 
S_{i}^{\mbox{\scriptsize{tot}}} \le \frac{\nu_{i} }{\pi^{2} V} .       
\end{equation} 

\noindent \textbf{Proof:} 
The proof of this Theorem in given in \citet{sobkuc09}. It is based on inequality:
\[\int_{0}^{1}u^{2} \left(x\right)dx \le \frac{1}{\pi^{2} } \int_{0}^{1}\left(\frac{\partial u}{\partial x} \right)^{2} dx \] 
and relationships \eqref{GrindEQ__15_} and \eqref{GrindEQ__17_}.
{\hfill\mbox{\rule{2 true mm}{3 true mm}}}

Consider the set of values $\nu_{1} ,...,\nu_{d} $, $1\le i  \le d$. One can expect that smaller $\nu_{i} $ correspond to less influential variables $x_{i} $. This importance criterion is similar to the modified Morris importance measure $\mu^{*} $, whose limiting values are 
\[\mu_{i}^{*} =\int_{H^{d} }\left|\frac{\partial G(\mb{x})}{\partial x_{i} } \right|d\mb{x} .     \] 

From a practical point of view the criteria $\mu_{i}$ and $\nu_{i}$ are equivalent: they are evaluated by the same numerical algorithm and are linked by relations $\nu_{i} \le C\mu_{i}$ and $\mu_{i} \le \sqrt{\nu_{i} } $. 

 The right term in \eqref{GrindEQ__32_} is further called the upper bound number one (UB1).

\noindent \textbf{Theorem 5. }There exists the following upper bound between DGSM \eqref{GrindEQ__12_} and the Sobol' total sensitivity index:
\begin{equation} \label{GrindEQ__33_} 
S_{i}^{\mbox{\scriptsize{tot}}} \le \frac{\varsigma_{i} }{V}.
\end{equation} 

\noindent \textbf{Proof:} The following inequality \citep{harlit73} is used:
\begin{equation} \label{GrindEQ__34_} 
0\le \int_{0}^{1}u^{2} dx -\left(\int_{0}^{1}udx \right)^{2} \le \frac{1}{2} \int_{0}^{1}x(1-x)u'^{2} dx .    
\end{equation} 
The inequality is reduced to an equality only if $u$ is constant. Assume that $u$ is given by \eqref{GrindEQ__13_}, then $\displaystyle \int_{0}^{1}udx =0$.
From \eqref{GrindEQ__34_}, equation \eqref{GrindEQ__33_} is obtained.
{\hfill\mbox{\rule{2 true mm}{3 true mm}}}

Further $\varsigma_{i} / D$ is called the upper bound number two (UB2). Note that $\frac{1}{2} x_{i} (1-x_{i} )$ for $0\le x_{i} \le 1$ is bounded: $0\le \frac{1}{2} x_{i} (1-x_{i} )\le 1/8 $. Therefore, $0\le \varsigma_{i} \le \nu_{i} / 8$.

\subsection{Normally distributed variables}


\subsubsection{Lower bound on $S_{i}^{\mbox{\scriptsize{tot}}}$}

\noindent \textbf{Theorem 6.} If $X_{i}$ is normally distributed with a mean $\mu_i$ and a finite variance $\sigma_{i}^{2}$, there exists the following lower bound between DGSM \eqref{GrindEQ__37_} and the Sobol' total sensitivity index:
\begin{equation} \label{GrindEQ__38_} 
\frac{\sigma_{i}^{4}}{(\mu_i^2+\sigma_i^2)V} w_{i}^{2} \le S_{i}^{\mbox{\scriptsize{tot}}} . 
\end{equation} 

\noindent \textbf{Proof:} 
Using the equation \eqref{GrindEQ__17_} and Cauchy-Schwartz inequality applied on $\displaystyle \int_{R^{d}} x_{i} u_{i} (\mb{x})d F(\mb{x})$ (with $F$ the joint cdf), \citet{kucson15} give the proof of this inequality when $\mu_i=0$ (omitting to mention this condition).
The general proof, obtained by \citet{pet15}, is given below.

Consider a univariate function $g(X)$, with $X$ a normally distributed variable with mean $\mu$, finite variance $\sigma^2$ and cdf $F$.
With adequate conditions on $g$, the following equality is obtained by integrating by parts:
\begin{eqnarray*}
& \mathbb{E}[g'(X)] = \displaystyle \int_{-\infty}^{\infty} g'(x) dF(x) = \displaystyle \frac{1}{\sigma \sqrt{2\pi}} \int_{-\infty}^{\infty} g'(x) \exp \left[-\frac{(x-\mu)^2}{2\sigma^2}\right] dx \\
& = \displaystyle \frac{1}{\sigma \sqrt{2\pi}} \left[ g(x) \exp \left[-\frac{(x-\mu)^2}{2\sigma^2} \right] \right]_{-\infty}^{+\infty} + \frac{1}{\sigma \sqrt{2\pi}} \int_{-\infty}^{\infty} g(x) \frac{x-\mu}{\sigma^2} \exp \left[-\frac{(x-\mu)^2}{2\sigma^2}\right] dx  \\
& = \displaystyle \frac{1}{\sigma^2} \int_{-\infty}^{\infty} x g(x) dF(x) - \mu \int_{-\infty}^{\infty} g(x) dF(x) .
\end{eqnarray*}

In this equation, replacing $g(x)$ by $u_i(\mb{x})$ with $x_i$ normally distributed, the $w_i$ DGSM  writes
\[ w_i =  \int_{R^{d}} \frac{\partial G(\mb{x})}{\partial x_{i}} dF(\mb{x}) =  \int_{R^{d}} \frac{\partial u_i(\mb{x})}{\partial x_{i}} dF(\mb{x}) = \frac{1}{\sigma^2_i} \int_{R^{d}} x_i u_i(\mb{x}) dF(\mb{x}), \]
because $\int_{R^{d}} u_i(\mb{x}) dF(\mb{x}) = 0$ (due to the ANOVA decomposition condition).
Moreover, the Cauchy-Schwartz inequality applied on $\int_{R^{d}} x_i u_i(\mb{x}) dF(\mb{x})$ gives
\[ \left[ \int_{R^{d}}  x_i u_i(\mb{x}) dF(\mb{x}) \right]^2 \le \int_{R^{d}} x_i^2 dF(\mb{x})  \int_{R^{d}} [u_i(\mb{x})]^2 dF(\mb{x}) .\]
Combining the two latter equations leads to the expression
\[ w_i^2 \le \frac{1}{\sigma^4_i}(\mu_i^2+\sigma_i^2) V S_{i}^{\mbox{\scriptsize{tot}}}, \]
which is equivalent to Eq. (\ref{GrindEQ__38_}).
{\hfill\mbox{\rule{2 true mm}{3 true mm}}}

%

\subsubsection{Upper bounds on $S_{i}^{\mbox{\scriptsize{tot}}} $}

The following Theorem 7 is a generalization of Theorem 1.

\noindent \textbf{Theorem 7.} If $X_{i}$ has a finite variance $\sigma_{i}^{2} $ and $\displaystyle c\le \left|\frac{\partial G}{\partial x_{i} } \right|\le C$, then
\begin{equation} \label{GrindEQ__41_} 
\frac{\sigma_{i}^{2} c^{2} }{V} \le S_{i}^{\mbox{\scriptsize{tot}}} \le \frac{\sigma_{i}^{2} C^{2} }{V}.       
\end{equation} 
The constant factor $\sigma_{i}^{2} $ cannot be improved.

\noindent \textbf{Theorem 8. }If $X_{i}$ is normally distributed with a finite variance $\sigma_{i}^{2} $, there exists the following upper bound between DGSM \eqref{GrindEQ__36_} and the Sobol' total sensitivity index:
\begin{equation} \label{GrindEQ__42_} 
S_{i}^{\mbox{\scriptsize{tot}}} \le \frac{\sigma_{i}^{2} }{V} \nu_{i} .       
\end{equation} 
The constant factor $\sigma_{i}^{2} $ cannot be reduced.

\noindent \textbf{Proof:} 
The proofs of these Theorems are presented in \citet{sobkuc09}.
{\hfill\mbox{\rule{2 true mm}{3 true mm}}}

\section{DGSM-based bounds for groups of variables}

Let $\mb{x}=\left(x_{1} ,...,x_{d} \right)$ be a point in the $d-$dimensional unit hypercube with Lebesgue measure $d\mb{x}=dx_{1} \cdot \cdot \cdot dx_{d} $. Consider an arbitrary subset of the variables $y=\left(x_{i_{1} } ,...,x_{i_{s} } \right)$, $1 \le i_1 \le \ldots \le i_s \le d$, and the set of remaining complementary variables $z$, so that $\mb{x}=(y,z)$, $d\mb{x}=dy \,dz$. 
Further all the integrals are written without integration limits, by assuming that each integration variable varies independently from $0$ to $1$.

Consider the following DGSM $\tau_{y} $:
\begin{equation} \label{GrindEQ__49_} 
\tau_{y} =\sum_{p=1}^{s}\int \left(\frac{\partial G\left(\mb{x}\right)}{\partial x_{i_{p} } } \right)^{2} \frac{1-3x_{i_{p} } +3x_{i_{p} }^{2} }{6}   d\mb{x}.     
\end{equation} 
\noindent \textbf{Theorem 9. }If $G\left(\mb{x}\right)$ is linear with respect to $x_{i_{1} } ,...,x_{i_{s} } $, then $V_{y}^{\mbox{\scriptsize{tot}}} =\tau_{y} $, or in other words $\displaystyle S_{y}^{\mbox{\scriptsize{tot}}} =\frac{\tau_{y} }{V} $.

\noindent \textbf{Theorem 10. } The following general inequality holds: $\displaystyle V_{y}^{\mbox{\scriptsize{tot}}} \le \left({24\mathord{\left/ {\vphantom {24 \pi^{2} }} \right. \kern-\nulldelimiterspace} \pi^{2} } \right)\tau_{y} $, or in other words $\displaystyle S_{y}^{\mbox{\scriptsize{tot}}} \le \frac{24}{\pi^{2}V} \tau_{y}$.

\noindent \textbf{Proof:} 
The proofs of these Theorems are given in \citet{sobkuc10}. The second theorem shows that small values of $\tau_{y} $ imply small values of $S_{y}^{\mbox{\scriptsize{tot}}} $ and this allows identification of a set of unessential factors $y$ (usually defined by a condition of the type $S_{y}^{\mbox{\scriptsize{tot}}} <\epsilon $, where $\epsilon $ is small).
{\hfill\mbox{\rule{2 true mm}{3 true mm}}}

\subsection{Importance criterion $\tau_{i} $}

Consider the one dimensional case when the subset $y$ consists of only one variable $y=\left(x_{i} \right)$, then measure $\tau_{y}=\tau_i $ has the form 
\begin{equation} \label{GrindEQ__50_} 
\tau_{i} =\int \left(\frac{\partial G\left(\mb{x}\right)}{\partial x_{i} } \right)^{2} \frac{1-3x_{i} +3x_{i}^{2} }{6} d\mb{x} .    
\end{equation} 
It is easy to show that $\nu_i/24 \leq \tau_i \leq \nu_i/6$. From UB1 it follows that 
\begin{equation} \label{GrindEQ__51_} 
S_{i}^{\mbox{\scriptsize{tot}}} \le \frac{24}{\pi^{2} V} \tau_{i} .        
\end{equation} 
Thus small values of $\tau_{i} $ imply small values of $S_{i}^{\mbox{\scriptsize{tot}}} $, that are characteristic for non important variables $x_{i} $. 
At the same time, the following corollary is obtained from Theorem 9: if $G\left(\mb{x}\right)$ depends linearly on $x_{i} $, then $S_{i}^{\mbox{\scriptsize{tot}}} =\tau_{i} / V $. Thus $\tau_{i} $ is closer to $V_{i}^{\mbox{\scriptsize{tot}}} $ than $\nu_{i} $.

Note that the constant factor $1/\pi^2$ in \eqref{GrindEQ__32_} is the best possible. But in the general inequality for $\tau_{i} $ \eqref{GrindEQ__51_} the best possible constant factor is unknown.

There is a general link between importance measures $\tau_{i} $, $\varsigma_{i}$ and  $\nu_{i} $: 
\[\tau_{i} =-\varsigma_{i} +\frac{1}{6} \nu_{i} ,      \] 
then 
\[\varsigma_{i} =\frac{1}{6} \nu_{i} -\tau_{i} .      \]

\subsection{Normally distributed random variables}

Consider independent normal random variables $X_{1} ,...,X_{d} $ with parameters $(\mu_{i},\sigma_{i} )_{i=1\ldots d}$. Define $\tau_{i} $ as
\[\tau_{i} =\frac{1}{2} \mathbb{E}\left[\left(\frac{\partial G\left(\mb{x}\right)}{\partial x_{i} } \right)^{2} \left(x'_{i} -x_{i} \right)^{2} \right].    \] 
The expectation over $x'_{i} $ can be computed analytically. Then 
\[\tau_{i} =\frac{1}{2} \mathbb{E}\left[\left(\frac{\partial G\left(\mb{x}\right)}{\partial x_{i} } \right)^{2} \frac{\left(x_{i} - \mu_{i} \right)^{2} + \sigma_{i}^{2} }{2} \right].   \] 

\noindent \textbf{Theorem 11.} If $X_{1} ,...,X_{d} $ are independent normal random variables, then for an arbitrary subset $y$ of these variables, the following inequality is obtained:
\[S_{y}^{\mbox{\scriptsize{tot}}} \le \frac{2}{V} \tau_{y}.      \] 

\noindent \textbf{Proof:} 
The proof is given in \citet{sobkuc10}.
{\hfill\mbox{\rule{2 true mm}{3 true mm}}}

\section{DGSM-based upper bounds in the general case}

As previously, consider the function $G\left(X_{1} ,...,X_{d} \right)$, where $X_{1} ,...,X_{d}$ are independent random variables, defined in the Euclidian space $R^{d}$, with cdfs $F_{1} \left(x_{1} \right),...,F_{d} \left(x_{d} \right)$.
Assume further that each $X_i$ admits a probability density function (pdf), denoted by $f_i(x_i)$.
In the following, all the integrals are written without integration limits.

The developments in this section are based on the classical $L^2$-Poincar\'e inequality:
\begin{equation}
\label{eq:Poincare}
\int G(\mb{x})^2 dF(\mb{x}) 
\leq C(F) \int  \Vert \nabla G(\mb{x}) \Vert ^2 dF(\mb{x})
\end{equation}
where $F$ is the joint cdf of $(X_{1} ,...,X_{d})$.
\eqref{eq:Poincare} is valid for all functions $G$ in $L^2(F)$ such that $\int G(\mb{x}) dF(\mb{x})=0$ and $\Vert \nabla f \Vert \in L^2(F)$.
The constant $C(F)$ in Eq. (\ref{eq:Poincare}) is called a Poincar\'e constant of $F$.
In some cases, it exists and optimal Poincar\'e constant $C_{\text{opt}}(F)$ which is the best possible constant.
In measure theory, the Poincar\'e constants are expressed as a function of so-called Cheeger constants \citep{bob99} which are used for SA in \citet{lamioo13} (see \citet{roufru14} for more details).

A connection between total indices and DGSM has been established by \citet{lamioo13} for variables with continuous distributions (called Boltzmann probability measures in their paper). 

\noindent \textbf{Theorem 12. }
Let $F_i$ and $f_i$ be respectively the cdf and the pdf of $X_i$, the following inequality is obtained:
\begin{equation}
\label{eq:InequalityTotal}
S_{i}^{\mbox{\scriptsize{tot}}} \leq \, \frac{C(F_i) }{V}\nu_{i},
\end{equation}
with $\nu_i$ the DGSM defined in Eq. (\ref{GrindEQ__36_}) and 
\begin{equation}
C(F_i) = 4 \left[ \underset{x \in \mathbb{R}}{\sup} \frac{\min \left(F_i(x), 1-F_i(x) \right) }{f_i(x)} \right]^2.
\end{equation}

\noindent \textbf{Proof:} 
This result comes from the direct application of the $L^2$-Poincar\'e inequality (\ref{eq:Poincare}) on $u_i(\mb{x})$ (see Eq. \eqref{GrindEQ__13_}).
{\hfill\mbox{\rule{2 true mm}{3 true mm}}}

In \citet{lamioo13} and \citet{roufru14}, the particular case of log-concave probability distribution has been developed.
It includes classical distributions as for instance the normal, exponential, Beta, Gamma and Gumbel distributions.
In this case, the constant writes
\begin{equation}
C(F_i) = \frac{1}{f_i(\tilde{m_i})^2}
\end{equation}
with $\tilde{m_i}$ the median of the distribution $F_i$.
This allows to obtain analytical expressions for $C(F_i)$ in several cases \citep{lamioo13}.
In the case of a log-concave truncated distribution on $[a,b]$, the constant writes \citep{roufru14}
\begin{equation}
\left(F_i(b) - F_i(a) \right)^2 / f_i \left( q_i \left( \frac{F_i(a)+F_i(b)}{2} \right) \right) ^2
\end{equation}
with $q_i(\cdot)$ the quantile function of $X_i$.
Table \ref{tab:poinconst} gives some examples of Poincar\'e constants for several well-known and often used probability distributions in practice.

 \begin{table}[!ht]
 \begin{center}
    \begin{tabular}{lcc}
\hline
Distribution & Poincar\'e constant & Optimal constant\\
\hline
Uniform  ${\mathcal U}[a \, b]$  & $(b-a)^2/\pi^2$ &  yes\\
Normal ${\mathcal N}(\mu, \sigma^2)$ & $\sigma^2$ & yes \\
Exponential ${\mathcal E}(\lambda)$, $\lambda>0$ & $\displaystyle \frac{4}{\lambda^2}$ & yes\\
Gumbel ${\mathcal G}(\mu,\beta)$, scale $\beta>0$ & $\displaystyle \left(\frac{2\beta}{\log 2}\right)^2$ & no \\
Weibull ${\mathcal W}(k, \lambda)$, shape $k \geq 1$, scale $\lambda > 0$\hspace{0.2cm} & $\displaystyle \left[\frac{2\lambda(\log 2)^{(1-k)/k}}{k}\right]^2$ & no \\
\hline
  \end{tabular}
    \caption{Poincar\'e constants for a few probability distributions.}
    \label{tab:poinconst}
 \end{center}
\end{table}

For studying second-order interactions, \citet{roufru14} have derived a similar to \eqref{eq:InequalityTotal} inequality based on the squared crossed derivatives of the function.
Assuming that second-order derivatives of $G$ are in $L^2(F)$, it uses the so-called crossed-DGSM
\begin{equation}
\label{eq:crossedDGSM}
 \nu_{ij} = \int \left(\frac{\partial^{2}G(\mb{x})}{\partial x_i\partial x_j}\right)^{2}dF(\mb{x}),
\end{equation}
introduced by \citet{fripop08}.
An inequality link is made with an extension of the total Sobol' sensitivity indices to general sets of variables (called superset importance or total interaction index) proposed by \citet{liuowe06}.
In the case of a pair of variables $\{X_i, X_j\}$, the superset importance is defined as
\begin{equation}
\label{eq:superset2}
V_{ij}^{\text{super}} = \sum_{I \supseteq \{i,j\}} V_I.
\end{equation}
The estimation methods of this total interaction index have also been studied by \citet{frurou14}.

\noindent \textbf{Theorem 13. }
For all pairs $\{i,j\}$ ($1 \leq i<j \leq d$),
\begin{equation}
\label{eq:TheInequality1}
V_{ij} \, \leq \, V_{ij}^{\textrm{super}} \leq \, C(F_i) C(F_j) \nu_{ij} .
\end{equation}
These inequalities with the corresponding Sobol' indices write
\begin{equation}
\label{eq:TheInequality2}
S_{ij} \, \leq \, S_{ij}^{\textrm{super}} \leq \,  \frac{C(F_i) C(F_j)}{V} \nu_{ij}.
\end{equation}

\citet{roufru14} have shown on several examples how to apply this result in order to detect pairs of inputs that do not interact together (see also \citet{muerou12} and \citet{frurou14} which use Sobol' indices).

\section{Computational costs}

All DGSM can be computed using the same set of partial derivatives $\displaystyle \frac{\partial G(\mb{x})}{\partial x_{i} } ,{\rm \; }i=1,...,d$. Evaluation of $\displaystyle \frac{\partial G(\mb{x})}{\partial x_{i} } $ can be done analytically for explicitly given easily-differentiable functions or numerically:
\begin{equation} \label{GrindEQ__35_} 
\frac{\partial G(\mb{x}^{*} )}{\partial x_{i} } =\frac{\left[G\left(x_{1}^{*} ,\ldots ,x_{i-1}^{*} ,x_{i}^{*} +\delta ,x_{i+1}^{*} ,\ldots ,x_{n}^{*} \right)-G\left(\mb{x}^{*} \right)\right]}{\delta } .
\end{equation} 
This is called a finite-difference scheme (see \uline{Variational Methods}) with $\delta$ which is a small increment.
There is a similarity with the elementary effect formula (2) of the Morris method which is however computed with large $\Delta$.

In the case of straightforward numerical estimations of all partial derivatives \eqref{GrindEQ__35_} and computation of integrals using MC or QMC methods, the number of required function evaluations for a set of all input variables is equal to $N(d+1)$, where $N$ is a number of sampled points. Computing LB1 also requires values of $G\left(0,\mb{z}\right),G\left(1,\mb{z}\right)$, while computing LB2 requires only values of $G\left(1,\mb{z}\right)$. In total, numerical computation of LB* for all input variables would require $N_{G}^{\mbox{\scriptsize{LB*}}} =N(d+1)+2Nd=N(3d+1)$ function evaluations. Computation of all upper bounds require $N_{G}^{\mbox{\scriptsize{UB}}} =N(d+1)$ function evaluations.
This is the same number that the number of function evaluations required for computation of $S_{i}^{\mbox{\scriptsize{tot}}} $ which is $N_{G}^{S} =N(d+1)$ \citep{salann09}.
 
However, the number of sampled points \textit{N} needed to achieve numerical convergence can be different for DGSM and $S_{i}^{\mbox{\scriptsize{tot}}} $. 
It is generally lower for the case of DGSM. 
Moreover, the numerical efficiency of the DGSM method can be significantly increased by using algorithmic differentiation in the adjoint (reverse) mode \citep{griwal08} (see also \uline{Variational Methods}). 
This approach allows estimating all derivatives at a cost independent of $d$, at most 4-6 times of that for evaluating the original function $G(\mb{x})$ \citep{janleo14}.

\section{Test cases}

In this section, three test cases are considered, in order to illustrate application of DGSM and their links with $S_{i}^{\mbox{\scriptsize{tot}}} $.

\noindent \textbf{Example 1.} Consider a linear with respect to $x_{i} $ function:
\[G(x)=a(\mb{z})x_{i} +b(\mb{z}).       \] 
For this function $S_{i} =S_{i}^{\mbox{\scriptsize{tot}}} $, $\displaystyle V_{i}^{\mbox{\scriptsize{tot}}} =\frac{1}{12} \int_{H^{d-1} }a^{2} (\mb{z})d\mb{z} $, $\displaystyle \nu_{i} =\int_{H^{d-1} }a^{2} (\mb{z})d\mb{z} $, $\displaystyle \mbox{LB1}=\frac{\left(\int_{H^{d} }\left(a^{2} (\mb{z})-2a^{2} (\mb{z})x_{i} \right)d\mb{z}dx_{i}  \right)^{2} }{4V\int_{H^{d-1} }a^{2} (\mb{z})d\mb{z} } =0$ and $\displaystyle \gamma (m)=\frac{(2m+1)m^{2} \left(\int_{H^{d-1} }a(\mb{z})d\mb{z} \right)^{2} }{4(m+2)^{2} (m+1)^{2} V} $. A maximum value of $\gamma (m)$ is attained at $m^{*} $=3.745, while $\displaystyle \gamma^{*} (m^{*} )=\frac{0.0401}{V} \left(\int a(\mb{z})d\mb{z} \right)^{2} $. The lower and upper bounds are $\displaystyle \mbox{LB*} \approx {\rm 0.48}S_{i}^{\mbox{\scriptsize{tot}}} $, $\displaystyle \mbox{UB1} \approx 1{\rm .22}S_{i}^{\mbox{\scriptsize{tot}}} $. $\displaystyle \mbox{UB2} =\frac{1}{12V} \int_{0}^{1}a(\mb{z})^{2} d\mb{z} =S_{i}^{\mbox{\scriptsize{tot}}} $.

\noindent For this test function UB2 $<$ UB1.

\noindent \textbf{Example 2.} Consider the so-called g-function which is often used in global SA for illustration purposes:
\[G(x)=\prod_{i=1}^{d}v_{i}  ,\] 
where $\displaystyle v_{i} =\frac{|4x_{i} -2|+a_{i} }{1+a_{i} } $, $a_{i} (i=1,...,d)$ are constants. It is easy to see that for this function $\displaystyle g_{i} (x_{i} )=(v_{i} -1)$, $u_{i} (x)=(v_{i} -1)\prod_{j=1,j\ne i}^{d}v_{j}  $ and as a result LB1=0. The total variance is $\displaystyle V=-1+\prod_{j=1}^{d}\left(1+\frac{1/3}{(1+a_{j} )^{2} } \right) $. The analytical values of $S_{i} $, $\displaystyle S_{i}^{\mbox{\scriptsize{tot}}} $ and LB2 are given in Table \ref{tab:1}.

\begin{table}
\caption{The analytical expressions for $S_{i} $, $S_{i}^{\mbox{\scriptsize{tot}}} $and LB2 for g-function.}\label{tab:1}
\centering
\begin{tabular}{|c|c|c|} \hline 
$S_{i} $ & $S_{i}^{\mbox{\scriptsize{tot}}} $ & $\gamma (m)$ \\ \hline 
$\displaystyle \frac{1/3}{(1+a_{i} )^{2} V} $ & $\displaystyle \frac{\frac{1/3}{(1+a_{i} )^{2} } \prod_{j=1,j\ne i}^{d}\left(1+\frac{1/3}{(1+a_{j} )^{2} } \right) }{V} $ & $\displaystyle \frac{(2m+1)\left[1-\frac{4\left(1-(1/2)^{m+1} \right)}{m+2} \right]^{2} }{(1+a_{i} )^{2} (m+1)^{2} V} $ \\ \hline 
\end{tabular}
\end{table}

By solving equation $\displaystyle \frac{d\gamma (m)}{dm} =0$, $m^{*} $=9.64 and $\displaystyle \gamma (m^{*} )=\frac{0.0772}{(1+a_{i} )^{2} V} $. It is interesting to note that $m^{*} $ does not depend on $a_{i} {\rm ,\; }i=1,2,...,d$ and \textit{d}. In the extreme cases: if $a_{i} \to \infty $ for all \textit{i}, $\displaystyle \frac{\gamma (m^{*} )}{S_{i}^{\mbox{\scriptsize{tot}}} } \to 0.257$, $\displaystyle \frac{S_{i} }{S_{i}^{\mbox{\scriptsize{tot}}} } \to 1$, while if $a_{i} \to 0$ for all \textit{i}, $\displaystyle \frac{\gamma (m^{*} )}{S_{i}^{\mbox{\scriptsize{tot}}} } \to \frac{0.257}{(4/3)^{d-1} } $, $\displaystyle \frac{S_{i} }{S_{i}^{\mbox{\scriptsize{tot}}} } \to \frac{1}{(4/3)^{d-1} } $. The analytical expression for $S_{i}^{\mbox{\scriptsize{tot}}} $, UB1 and UB2 are given in Table \ref{tab:2}.

\begin{table}
\caption{The analytical expressions for $S_{i}^{\mbox{\scriptsize{tot}}} $, UB1 and UB2 for g-function.}\label{tab:2}
\centering
\begin{tabular}{|c|c|c|} \hline 
$S_{i}^{\mbox{\scriptsize{tot}}} $ & UB1 & UB2 \\ \hline 
$\displaystyle \frac{\frac{1/3}{(1+a_{i} )^{2} } \prod_{j=1,j\ne i}^{d}\left(1+\frac{1/3}{(1+a_{j} )^{2} } \right) }{V} $ & $\displaystyle \frac{16\prod_{j=1,j\ne i}^{d}\left(1+\frac{1/3}{(1+a_{j} )^{2} } \right) }{(1+a_{i} )^{2} \pi^{2} V} $ & $\displaystyle \frac{4\prod_{j=1,j\ne i}^{d}\left(1+\frac{1/3}{(1+a_{j} )^{2} } \right) }{3(1+a_{i} )^{2} V} $ \\ \hline 
\end{tabular}
\end{table}

For this test function $\displaystyle \frac{S_{i}^{\mbox{\scriptsize{tot}}} }{\mbox{UB1}} =\frac{\pi^{2} }{48} $, $\displaystyle \frac{S_{i}^{\mbox{\scriptsize{tot}}} }{\mbox{UB2}} =\frac{1}{4} $, hence $\displaystyle \frac{\mbox{UB2}}{\mbox{UB1}} =\frac{\pi^{2} }{12} <1$.

\begin{figure}
\centering
\includegraphics*[width=3.16in, height=2.50in, keepaspectratio=false]{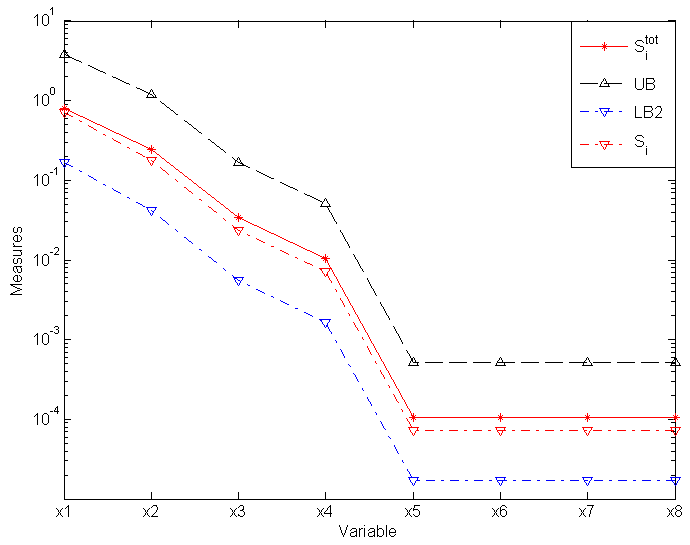}
\caption{Values of $S_{i} $, $S_{i}^{\mbox{\scriptsize{tot}}}$, LB2 and UB1 for all input variables. Example 2 with $\mb{a}=[0,1,4.5,9,99,99,99,99]$, $d=8$.}\label{fig:1}
\end{figure}

Values of $S_{i} $, $S_{i}^{\mbox{\scriptsize{tot}}} $, UB1, UB2 and LB2 for the case of \textbf{\textit{a}}=[0,1,4.5,9,99,99,99,99], \textit{d}=8 are given in Table \ref{tab:3} and shown in Fig. \ref{fig:1}. One can see that the knowledge of LB2 and UB1 allows to rank correctly all the variables in the order of their importance.

\begin{table}
\caption{Values of LB*, $S_{i} $, $S_{i}^{\mbox{\scriptsize{tot}}} $, UB1 and UB1. Example 2, \textbf{\textit{a}}=[0,1,4.5,9,99,99,99,99], \textit{d}=8.}\label{tab:3}
\centering
\begin{tabular}{|c|c|c|c|c|c|} \hline 
 & $x_{1}$ & $x_{2}$ & $x_{3}$ & $x_{4}$ & $x_{5} ...x_{8}$ \\ \hline 
LB* & $0.166$ & $0.0416$ & $0.00549$ & $0.00166$ & $0.000017$ \\ \hline 
$S_{i}$ & $0.716$ & $0.179$ & $0.0237$ & $0.00720$ & $0.0000716$ \\ \hline 
$S_{i}^{\mbox{\scriptsize{tot}}}$ & $0.788$ & $0.242$ & $0.0343$ & $0.0105$ & $0.000105$ \\ \hline 
UB1 & $3.828$ & $1.178$ & $0.167$ & $0.0509$ & $0.00051$ \\ \hline 
UB2 & $3.149$ & $0.969$ & $0.137$ & $0.0418$ & $0.00042$ \\ \hline 
\end{tabular}
\end{table}

\noindent \textbf{Example 3.} Consider the reduced Morris' test function with four inputs \citep{camcar07}:
\begin{equation}
	f(\boldsymbol{x})=\sum_{i=1}^4b_ix_i+\sum_{i\leq j}^4b_{ij}x_ix_j+\sum_{i\leq j\leq k}^4b_{ijk}x_ix_jx_k
\end{equation}
\begin{equation*}
\mbox{with }	b_i=\left[
	\begin{array}[c]{c}
	0.05 \\ 
	0.59 \\ 
	10 \\ 
	0.21
	\end{array}
	\right]
	\hspace{0.5cm},\hspace{0.5cm}
	b_{ij}=\left[
						\begin{array}[c]{c c c c}
							0&80&60&40\\ 
							0&30&0.73&0.18\\ 
							0&0&0.64&0.93\\ 
							0&0&0&0.06
						\end{array}
		\right]
	\hspace{0.5cm},\hspace{0.5cm}
	b_{ij4}=\left[
						\begin{array}[c]{c c c c}
							0&10&0.98&0.19\\ 
							0&0&0.49&50\\ 
							0&0&0&1\\ 
							0&0&0&0
						\end{array}
		\right] \;.
\end{equation*}
The indices $b_{ijk}\ \forall\ k\not=4$ are null. 

\begin{figure}[!ht]
$$  \includegraphics[width=0.7\textwidth]{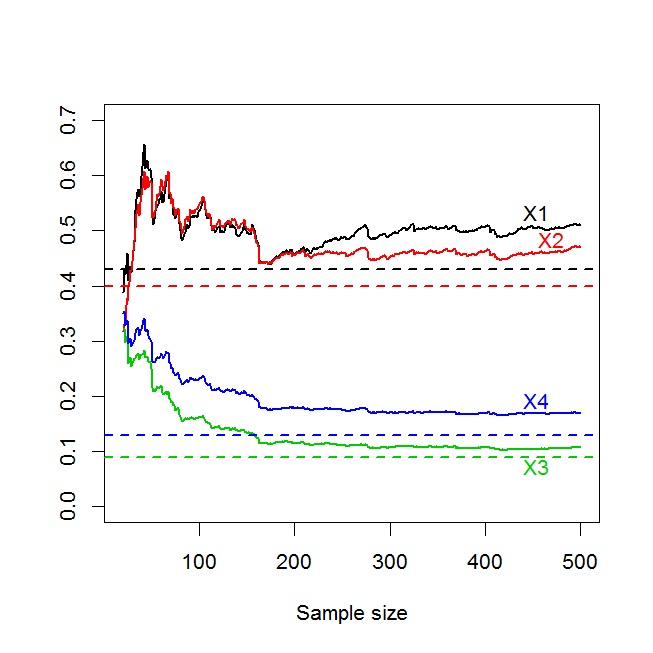}$$
\caption{For the $4$ input variables of the reduced Morris' test function: Convergence of the DGSM bound estimates (solid lines) in function of the sample size and comparison to theoretical values of total Sobol' indices $S_{T_i}$ (dashed lines).}
  \label{majoration_Morris}
 \end{figure}

The four input variables $x_i$ $(i=1,\ldots,4)$ follow uniform distribution on $[0,1]$. 
Sobol' indices are computed via the Monte-carlo scheme of \citet{sal02} (using two initial matrices of size $10^5$), while DGSM are computed with Monte-Carlo sampling of size $n$ (using derivatives computing by finite differences \eqref{GrindEQ__35_} with $\delta=10^{-5}$), with $n$ ranging from $20$ to $500$, Figure ~\ref{majoration_Morris} shows that DGSM bounds UB$1_i$ are greater than the total Sobol' indices $S_{T_i}$ (for $i=1,2,3,4$) as expected, except for $n<30$ which is a too small sample size.
For small $S_{T_i}$, UB$1_i$ is close to the $S_{T_i}$ value.
It confirms that DGSM bounds are first useful for screening exercises.
Other numerical tests involving non-uniform and non-normal distributions for the inputs can be found in \citet{lamioo13} and \citet{frurou14}.


\section{Conclusions}

This paper has shown that using lower and upper bounds based on DGSM is possible in most cases to get a good practical estimation of the values of $S_{i}^{\mbox{\scriptsize{tot}}}$ at a fraction of the CPU cost for estimating $S_{i}^{\mbox{\scriptsize{tot}}}$. 
Upper and lower bounds can be estimated using MC/QMC integration methods using the same set of partial derivative values.
Most of the applications show that DGSM can be used for fixing unimportant variables and subsequent model reduction because small values of DGSM imply small values of $S_{i}^{\mbox{\scriptsize{tot}}}$.
In a general case variable ranking can be different for DGSM and variance based methods but for linear function and product function, DGSM can give the same variable ranking as $S_{i}^{\mbox{\scriptsize{tot}}}$. 

Engineering applications of DGSM can be found for instance in \citet{kipkuc09} and \citet{rodban12} for biological systems modeling, \citet{patpra10b} for structural mechanics, \citet{ioopop12} for an aquatic prey-predator model, \citet{pet15} for a river flood model and \citet{toubus14} for an hydrogeological simulator of the oil industry.
One of the main prospect in practical situations is to use algorithmic differentiation in the reverse (adjoint) mode on the numerical model, allowing to estimate efficiency all partial derivatives of this model (see \uline{Variational Methods}). 
In this case, the cost of DGSM estimations would be independent of the number of input variables.
Obtaining global sensitivity information in a reasonable cpu time cost is therefore possible even for large-dimensional model (several tens and spatially distributed inputs in the recent and pioneering attempt of \citet{pet15}).
When the adjoint model is not available, the DGSM estimation remains a problem in high dimension and novel ideas have to be explored \citep{patpra10a} \citep{patpra10b}.
Coupling DGSM with non-parametric regression techniques or metamodel-based technique (see \uline{Metamodel-based sensitivity analysis: Polynomial chaos expansions and Gaussian processes}) is another research prospect as first shown by \citet{sudmai15} and \citet{delmar15}.


The authors would like to thank Prof. I. Sobol', Dr. S. Song, S. Petit, Dr. M. Lamboni, Dr. O. Roustant and Prof. F. Gamboa for their contributions to this work.
One of the authors (SK) gratefully acknowledges the financial support by the EPSRC grant EP/H03126X/1.


\bibliographystyle{spbasic}

\begin{thebibliography}{0}
\providecommand{\natexlab}[1]{#1}
\providecommand{\url}[1]{{#1}}
\providecommand{\urlprefix}{URL }
\expandafter\ifx\csname urlstyle\endcsname\relax
  \providecommand{\doi}[1]{DOI~\discretionary{}{}{}#1}\else
  \providecommand{\doi}{DOI~\discretionary{}{}{}\begingroup
  \urlstyle{rm}\Url}\fi
\providecommand{\eprint}[2][]{\url{#2}}

\end{thebibliography}


\begin{thebibliography}{41}
\providecommand{\natexlab}[1]{#1}
\providecommand{\url}[1]{{#1}}
\providecommand{\urlprefix}{URL }
\expandafter\ifx\csname urlstyle\endcsname\relax
  \providecommand{\doi}[1]{DOI~\discretionary{}{}{}#1}\else
  \providecommand{\doi}{DOI~\discretionary{}{}{}\begingroup
  \urlstyle{rm}\Url}\fi
\providecommand{\eprint}[2][]{\url{#2}}

\bibitem[{Bobkov(1999)}]{bob99}
Bobkov SG (1999) Isoperimetric and analytic inequalities for log-concave
  probability measures. The Annals of Probability 27(4):1903--1921

\bibitem[{Campolongo and Braddock(1999)}]{cambra99}
Campolongo F, Braddock R (1999) The use of graph theory in the sensitivity
  analysis of model output: a second order screening method. Reliability
  Engineering and System Safety 64:1--12

\bibitem[{Campolongo et~al(2007)Campolongo, Cariboni, and Saltelli}]{camcar07}
Campolongo F, Cariboni J, Saltelli A (2007) An effective screening design for
  sensitivity analysis of large models. Environmental Modelling and Software
  22:1509--1518

\bibitem[{Cropp and Braddock(2002)}]{crobra02}
Cropp R, Braddock R (2002) The new {M}orris method: an efficient second-order
  screening method. Reliability Engineering and System Safety 78:77--83

\bibitem[{{De Lozzo} and Marrel(2015)}]{delmar15}
{De Lozzo} M, Marrel A (2015) Estimation of the derivative-based global
  sensitivity measures using a {G}aussian process metamodel. Submitted

\bibitem[{F\'edou and Rendas(2015)}]{fedran15}
F\'edou JM, Rendas MJ (2015) Extending {M}orris method: identification of the
  interaction graph using cycle-equitable designs. Journal of Statistical
  Computation and Simulation 85:1398--1419

\bibitem[{Friedman and Popescu(2008)}]{fripop08}
Friedman J, Popescu B (2008) {Predictive Learning via Rule Ensembles}. The
  Annals of Applied Statistics 2(3):916--954

\bibitem[{Fruth et~al(2014)Fruth, Roustant, and Kuhnt}]{frurou14}
Fruth J, Roustant O, Kuhnt S (2014) Total interaction index: {A} variance-based
  sensitivity index for second-order interaction screening. Journal of
  Statistical Planning and Inference 147:212--223

\bibitem[{Griewank and Walther(2008)}]{griwal08}
Griewank A, Walther A (2008) Evaluating derivatives: {P}rinciples and
  techniques of automatic differentiation. SIAM Philadelphia

\bibitem[{Hardy et~al(1973)Hardy, Littlewood, and Polya}]{harlit73}
Hardy G, Littlewood J, Polya G (1973) Inequalities. Cambridge University Press,
  Second edition

\bibitem[{Homma and Saltelli(1996)}]{homsal96}
Homma T, Saltelli A (1996) Importance measures in global sensitivity analysis
  of non linear models. Reliability Engineering and System Safety 52:1--17

\bibitem[{Iooss et~al(2012)Iooss, Popelin, Blatman, Ciric, Gamboa, Lacaze, and
  Lamboni}]{ioopop12}
Iooss B, Popelin AL, Blatman G, Ciric C, Gamboa F, Lacaze S, Lamboni M (2012)
  Some new insights in derivative-based global sensitivity measures. In:
  Proceedings of the PSAM11 ESREL 2012 Conference, Helsinki, Finland, pp
  1094--1104

\bibitem[{Jansen et~al(2014)Jansen, Leovey, Nube, Griewank, and
  {Mueller-Preussker}}]{janleo14}
Jansen K, Leovey H, Nube A, Griewank A, {Mueller-Preussker} M (2014) A first
  look of quasi-{M}onte {C}arlo for lattice field theory problems. Computer
  Physics Communication 185:948--959

\bibitem[{Jansen(1999)}]{jan99}
Jansen M (1999) Analysis of variance designs for model output. Computer Physics
  Communication 117:25--43

\bibitem[{Kiparissides et~al(2009)Kiparissides, Kucherenko, Mantalaris, and
  Pistikopoulos}]{kipkuc09}
Kiparissides A, Kucherenko S, Mantalaris A, Pistikopoulos E (2009) Global
  sensitivity analysis challenges in biological systems modeling. Journal of
  Industrial and Engineering Chemistry Research 48:1135--1148

\bibitem[{Kucherenko and Song(2015)}]{kucson15}
Kucherenko S, Song S (2015) Derivative-based global sensitivity measures and
  their link with {S}obol' sensitivity indices. In: Cools R, Nuyens D (eds)
  Proceedings of the Eleventh International Conference on Monte Carlo and
  Quasi-Monte Carlo Methods in Scientific Computing (MCQMC 2014),
  Springer-Verlag,, Leuven, Belgium

\bibitem[{Kucherenko et~al(2009)Kucherenko, Rodriguez-Fernandez, Pantelides,
  and Shah}]{kucrod09}
Kucherenko S, Rodriguez-Fernandez M, Pantelides C, Shah N (2009) Monte carlo
  evaluation of derivative-based global sensitivity measures. Reliability
  Engineering and System Safety 94:1135--1148

\bibitem[{Lamboni(2013)}]{lam13}
Lamboni M (2013) New way of estimating total sensitivity indices. In:
  Proceedings of the 7th International Conference on Sensitivity Analysis of
  Model Output (SAMO 2013), Nice, France

\bibitem[{Lamboni et~al(2013)Lamboni, Iooss, Popelin, and Gamboa}]{lamioo13}
Lamboni M, Iooss B, Popelin AL, Gamboa F (2013) Derivative-based global
  sensitivity measures: general links with sobol' indices and numerical tests.
  Mathematics and Computers in Simulation 87:45--54

\bibitem[{Liu and Owen(2006)}]{liuowe06}
Liu R, Owen A (2006) Estimating mean dimensionality of analysis of variance
  decompositions. Journal of the American Statistical Association
  101(474):712--721

\bibitem[{Morris(1991)}]{mor91}
Morris M (1991) Factorial sampling plans for preliminary computational
  experiments. Technometrics 33:161--174

\bibitem[{Muehlenstaedt et~al(2012)Muehlenstaedt, Roustant, Carraro, and
  Kuhnt}]{muerou12}
Muehlenstaedt T, Roustant O, Carraro L, Kuhnt S (2012) Data-driven {K}riging
  models based on {FANOVA}-decomposition. {Statistics \& Computing} 22:723--738

\bibitem[{Patelli and Pradlwarter(2010)}]{patpra10a}
Patelli E, Pradlwarter H (2010) Monte {C}arlo gradient estimation in high
  dimensions. International Journal for Numerical Methods in Engineering
  81:172--188

\bibitem[{Patelli et~al(2010)Patelli, Pradlwarter, and Schu\"eller}]{patpra10b}
Patelli E, Pradlwarter HJ, Schu\"eller GI (2010) Global sensitivity of
  structural variability by random sampling. Computer Physics Communications
  181:2072--2081

\bibitem[{Petit(2015)}]{pet15}
Petit S (2015) Analyse de sensibilit\'e globale du module {MASCARET} par
  l'utilisation de la diff\'erentiation automatique. Rapport de stage de fin
  d'\'etudes de Sup\'elec, EDF R\&D, Chatou, France

\bibitem[{Pujol(2009)}]{puj09}
Pujol G (2009) Simplex-based screening designs for estimating metamodels.
  Reliability Engineering and System Safety 94:1156--1160

\bibitem[{Rodriguez-Fernandez et~al(2012)Rodriguez-Fernandez, Banga, and
  Doyle}]{rodban12}
Rodriguez-Fernandez M, Banga J, Doyle F (2012) Novel global sensitivity
  analysis methodology accounting for the crucial role of the distribution of
  input parameters: application to systems biology models. International
  Journal of Robust Nonlinear Control 22:1082--1102

\bibitem[{Roustant et~al(2014)Roustant, Fruth, Iooss, and Kuhnt}]{roufru14}
Roustant O, Fruth J, Iooss B, Kuhnt S (2014) Crossed-derivative-based
  sensitivity measures for interaction screening. Mathematics and Computers in
  Simulation 105:105--118

\bibitem[{Saltelli(2002)}]{sal02}
Saltelli A (2002) Making best use of model evaluations to compute sensitivity
  indices. Computer Physics Communication 145:280--297

\bibitem[{Saltelli et~al(2008)Saltelli, Ratto, Andres, Campolongo, Cariboni,
  Gatelli, Salsana, and Tarantola}]{salrat08}
Saltelli A, Ratto M, Andres T, Campolongo F, Cariboni J, Gatelli D, Salsana M,
  Tarantola S (2008) Global sensitivity analysis - The primer. Wiley

\bibitem[{Saltelli et~al(2010)Saltelli, Annoni, Azzini, Campolongo, Ratto, and
  Tarantola}]{salann09}
Saltelli A, Annoni P, Azzini I, Campolongo F, Ratto M, Tarantola S (2010)
  Variance based sensitivity analysis of model output. {D}esign and estimator
  for the total sensitivity index. Computer Physics Communication 181:259--270

\bibitem[{Santiago et~al(2012)Santiago, Corre, Claeys-Bruno, and
  Sergent}]{sancor12}
Santiago J, Corre B, Claeys-Bruno M, Sergent M (2012) Improved sensitivity
  through {M}orris extension. Chemometrics and Intelligent Laboratory Systems
  113:52--57

\bibitem[{Sobol(1990)}]{sob90}
Sobol I (1990) Sensitivity estimates for non linear mathematical models (in
  {R}ussian). Matematicheskoe Modelirovanie 2:112--118

\bibitem[{Sobol(1993)}]{sob93}
Sobol I (1993) Sensitivity estimates for non linear mathematical models.
  Mathematical Modelling and Computational Experiments 1:407--414

\bibitem[{Sobol(2001)}]{sob01}
Sobol I (2001) Global sensitivity indices for non linear mathematical models
  and their {M}onte {C}arlo estimates. Mathematics and Computers in Simulation
  55:271--280

\bibitem[{Sobol and Gershman(1995)}]{sobger95}
Sobol I, Gershman A (1995) On an alternative global sensitivity estimators. In:
  Proceedings of SAMO 1995, Belgirate, pp 40--42

\bibitem[{Sobol and Kucherenko(2005)}]{sobkuc05}
Sobol I, Kucherenko S (2005) Global sensitivity indices for non linear
  mathematical models. {R}eview. Wilmott Magazine 1:56--61

\bibitem[{Sobol and Kucherenko(2009)}]{sobkuc09}
Sobol I, Kucherenko S (2009) Derivative based global sensitivity measures and
  their links with global sensitivity indices. Mathematics and Computers in
  Simulation 79:3009--3017

\bibitem[{Sobol and Kucherenko(2010)}]{sobkuc10}
Sobol I, Kucherenko S (2010) A new derivative based importance criterion for
  groups of variables and its link with the global sensitivity indices.
  Computer Physics Communications 181:1212 -- 1217

\bibitem[{Sudret and Mai(2015)}]{sudmai15}
Sudret B, Mai CV (2015) Computing derivative-based global sensitivity measures
  using polynomial chaos expansions. Reliability Engineering and System Safety
  134:241--250

\bibitem[{Touzany and Busby(2014)}]{toubus14}
Touzany S, Busby D (2014) Screening method using the derivative-based global
  sensitivity indices with application to reservoir simulator. Oil \& Gas
  Science and Technology – Rev IFP Energies nouvelles 69:619--632

\end{thebibliography}

\end{document}